\documentclass[sn-mathphys,Numbered]{sn-jnl}% Math and Physical Sciences Reference Style
\usepackage{comment}
\usepackage{graphicx}%
\usepackage{multirow}%
\usepackage{amsmath,amssymb,amsfonts}%
\usepackage{amsthm}%
\usepackage{mathrsfs}%
\usepackage[title]{appendix}%
\usepackage{xcolor}%
\usepackage{textcomp}%
\usepackage{manyfoot}%
\usepackage{booktabs}%
\usepackage{algorithm}%
\usepackage{algorithmicx}%
\usepackage{algpseudocode}%
\usepackage{listings}%
\usepackage{mathtools}
\usepackage{amsmath}
\usepackage{amsfonts}
\usepackage{amsthm}
\usepackage{multicol,multirow}
\usepackage{nicefrac}
\usepackage{color, colortbl}
\definecolor{green}{rgb}{0.376,0.796,0.517}

\def \PP {\textbf P}
\def \EE {\mathbb E}
\def \R {\mathbb R}

\def\e{\varepsilon}

\newtheorem{assumption}{Assumption}
\DeclareMathOperator*{\argmax}{arg\,max}
\DeclareMathOperator*{\argmin}{arg\,min}
\DeclareMathOperator*{\esssup}{ess\,sup}

\newtheorem{theorem}{Theorem}%  meant for continuous numbers
\newtheorem{remark}{Remark}%
\newtheorem{lemma}{Lemma}%
\newtheorem{definition}{Definition}%

\begin{document}

\title[Fast convergence of sample-average approximation for saddle-point problems]{Fast convergence of sample-average approximation for saddle-point problems}


\author*[1]{\fnm{Vitali} \sur{Pirau}}\email{pireyvitalik@phystech.edu}

\affil*[1]{\orgdiv{Department of discrete mathematics}, \orgname{MIPT}, \orgaddress{\street{Institutskiy per.}, \city{Dolgoprudny}, \postcode{141701}, \state{Moscow region}, \country{Russia}}}


\abstract{Stochastic saddle point (SSP) problems are, in general, less studied compared to stochastic minimization problems. However, SSP problems emerge from machine learning (adversarial training, e.g., GAN, AUC maximization), statistics (robust estimation), and game theory (stochastic matrix games). Notwithstanding the existing results on the convergence of stochastic gradient algorithms, there is little analysis on the generalization, i.e.,  how learned on train data models behave with new data samples. In contrast to existing $\mathcal{O}(\frac{1}{n})$ in expectation result \cite{zhang2021generalization}, \cite{farnia2021train} and $\Tilde{\mathcal{O}}(\frac{\log(\nicefrac{1}{\delta})}{\sqrt{n}})$ results with high probability \cite{kang2022stability}, \cite{lei2021stability}, \cite{ji2021understanding}, we present $\mathcal{O}(\frac{\log(\nicefrac{1}{\delta})}{n})$ result in high probability for strongly convex-strongly concave setting. The main idea is local norms analysis. We illustrate our results in a matrix-game problem.}




\maketitle

\section{Introduction}
In practice, we are often faced with optimization problems, where the variables are naturally divided into two non-overlapping groups, one of which should be minimized, and the other maximized. This type of problem emerges from adversarial learning \cite{goodfellow2014generative}, reinforcement learning \cite{du2017stochastic}, \cite{dai2018sbeed}, AUC maximization \cite{zhao2011online}, \cite{gao2013one}, \cite{ying2016stochastic}, \cite{liu2018fast}, \cite{lei2021stability}, game theory \cite{facchinei2014non}, \cite{bach2019universal}. 

We will consider the minimax problem in the form:
\begin{equation}\label{ssp}
    \min\limits_{x \in \mathcal{X}}\max\limits_{y \in \mathcal{Y}} \left[F(x, y) = \PP\left(F(x, y, \xi)\right)\right],
\end{equation}
where $\PP F(x, y, \xi)$ denotes expectation of $F(x, y, \xi)$ with respect to $\xi \sim P$ and $\mathcal{X}, \mathcal{Y}$ are convex compact subsets of $\mathcal{R}^d$. The sample distribution $P$ is unknown.

The stochastic Approximation (SA) approach to solving the SSP problem is well-studied in the literature \cite{nemirovski2009robust}, \cite{zhao2022accelerated}. The primary focus of our work is to study both the Sample Average Approximation (SAA) approach \cite{shapiro2021lectures} to solve SSP and to give a high probability bound for the generalization of empirical solution for SSP. To our best knowledge, the SAA approach for SSP was studied only by \cite{zhang2021generalization}, where an expectation bound is provided. The generalization analysis is also provided in \cite{mehta2017fast}, \cite{zhang2017empirical}, \cite{wang2022stability}, \cite{farnia2021train}.

The empirical counterpart to \eqref{ssp} is:
\begin{equation}\label{esp}
    \min\limits_{x \in \mathcal{X}}\max\limits_{y \in \mathcal{Y}} \left[\hat{F}(x, y) = \PP_n\left(F(x, y, \xi)\right)\right],
\end{equation}
where $n$ is the number of data samples, $\PP_nF(x, y, \xi) = \frac{1}{n}\sum\limits_{i = 1}^n F(x, y, \xi_i)$, $\xi_i \sim P$. Denote the sample by $S_n = \{\xi_1, ..., \xi_n\}$.

We will further call this problem the Empirical Saddle Point (ESP) problem.

Where exists different measures of the quality of the SSP solutions. 

\begin{definition}
    Weak excess risk is defined as
    \begin{equation}\label{weakrisk}
        \Delta^w(x, y) = \sup_{y' \in \mathcal{Y}}\EE_{S_n}\left[F(x, y')\right] - \inf_{x' \in \mathcal{X}}\EE_{S_n}\left[F(x', y, \xi)\right].
    \end{equation}
\end{definition}

\begin{definition}
    Strong excess risk is defined as
    \begin{equation}\label{strongrisk}
        \Delta^s(x, y) = F(x, y^*(x)) - F(x^*(y), y).
    \end{equation}
\end{definition}

We will be interested in the high-probability bound of the strong excess risk \eqref{strongrisk}.

\textbf{Notation}. Throughout the paper, the notation $f_n = \mathcal{O}(g_n)$ means that for some constant $c$, which \textit{doesn't} depend on sample size $n$ and samples $\xi_i$, we have $f_n \leq c g_n$ for all natural $n$. We will denote a ball of radius $r$ with center in $x$ by $B_r(x)$. The space to which the ball belongs will be clear from the context.

\textbf{Contributions.} Our main result is high-probability generalization bound $\mathcal{O}(\frac{\log(\nicefrac{1}{\delta})}{n})$ for the ESP solution under strong convexity assumption. While \cite{zhang2021generalization} use stability argument and prove only in-expectation result (which can be continued with Markov inequality to $\mathcal{O}(\nicefrac{1}{n\delta})$ high-probability convergence), we explore local properties of the risk and get better result. Our analysis advances one described in \cite{lei2021stability} and gives a sharper bound, in contrast to previously best achieved $\mathcal{O}(\nicefrac{\sqrt{\log(n)}\log^2(\nicefrac{1}{\delta})}{\sqrt{n}\sigma})$ \cite{lei2021stability}.  We generalize local norms analysis from \cite{puchkin2023exploring} to minimax problems \eqref{ssp} and achieve similar results. Unfortunately, the described technique admits only a strongly convex-strongly concave setup, compared to the exp-concave setup, described in \cite{puchkin2023exploring} and requires bounded non-separated functional part (Assumption \eqref{lipgradbound}), which leads to the necessity of the tight regularization. The review of the recent results on generalization bounds for SSP problems and our contribution is summarized in Table \ref{tableresult}.

\begin{table}[h]
\caption{Summary of the results on stability of SSP}\label{tableresult}

\begin{tabular*}{\textwidth}{@{\extracolsep\fill}lcccccc}
\toprule
    \hline
     Algorithm  & Reference  &  Assumption & Measure & Rate \\\hline
     ERM & \cite{zhang2021generalization}& \ref{scsc}, \ref{lip} \ref{lipgrad} & Strong risk \eqref{strongrisk} &  $\mathcal{O}(\nicefrac{1}{n\sigma})$\\\hline
     SGDA & \cite{lei2021stability}& \ref{scsc}, \ref{lip} & Weak risk \eqref{weakrisk} &  $\mathcal{O}(\nicefrac{\sqrt{\log(n)}}{n\sigma})$\\\hline
     SGDA & \cite{lei2021stability}& C-$\sigma$-SC, \ref{lip} & High-probability strong risk \eqref{strongrisk} &  $\mathcal{O}(\nicefrac{\sqrt{\log(n)}\log^2(\nicefrac{1}{\delta})}{\sqrt{n}\sigma})$\\\hline
     SGDA & \cite{farnia2021train}& \ref{scsc}, \ref{lip}, \ref{lipgrad} & Weak risk \eqref{weakrisk} &  $\mathcal{O}(\nicefrac{\sqrt{\log(n)}}{\sqrt{n}\sigma})$\\\hline
     MC-SGDA & \cite{lei2021stability}& C-$\sigma$-SC, \ref{lip}, \ref{lipgrad} & Weak risk \eqref{weakrisk} &  $\mathcal{O}(\nicefrac{\sqrt{\log(n)}}{\sqrt{n}\sigma})$\\\hline
     \rowcolor{green}
     ERM & This work & $\sigma$-SC-SC, \ref{lip}, \ref{lipgrad}, \ref{lipgradbound} & High-probability strong risk \eqref{strongrisk} &  $\mathcal{O}(\nicefrac{\log(\nicefrac{1}{\delta})}{n\sigma})$\\\hline
\end{tabular*}
\footnotetext{Summary of Results. Bounds are stated in expectation or with high probability. For online setup, the optimal $T$ (number of iterations) is chosen to trade off generalization and optimization.  Here, C-C means convex-concave, C-$\sigma$-SC means convex-$\sigma$-strongly-concave. SGDA means Stochastic Gradient Descent Ascent.}
\end{table}

\section{Related work}

The SA approach for SSP problems is well-studied for convex \cite{nemirovski2009robust}, \cite{bach2019universal}, strongly convex \cite{natole2018stochastic}, \cite{yan2020optimal} and special (such as finite-sum and bilinear coupling) \cite{du2019linear}, \cite{shalev2013stochastic}, \cite{zhang2015stochastic} problems.

The pioneering work on stability and generalization is \cite{bousquet2002stability}, where the concept of uniform stability was introduced. \cite{shalev2010learnability} studied the concept of stability and showed its relation to the fundamental problem of learnability. The stability of gradient algorithms was further studied in multiple papers, e.g., \cite{hardt2016train}, \cite{lin2016generalization}, \cite{madden2020high}. For exp-concave functions stability of ERM was studied in \cite{mehta2017fast}, \cite{koren2015fast}, \cite{puchkin2023exploring}.

There is relatively little analysis of stability and generalization for saddle point problems. \cite{zhang2021generalization} firstly generalized  uniform stability argument \cite{bousquet2002stability} for empirical saddle point problems. The stability of SA algorithms for SSP was studied in \cite{zhang2017empirical}, \cite{wang2022stability}, \cite{farnia2021train}. In \cite{zhang2022uniform} generalization bounds for non-convex problems are provided.

\section{Assumptions and main result}
Let us denote by $(x^*, y^*)$ the solution to \eqref{ssp} and $x^*(y) = \argmin\limits_{x \in \mathcal{X}}F(x, y)$, $y^*(x) = \argmax\limits_{y \in \mathcal{Y}}F(x, y)$. Similarly, we will denote by $(\hat{x}, \hat{y})$ the solution to \eqref{esp} and $\hat{x}(y) = \argmin\limits_{x \in \mathcal{X}}\hat{F}(x, y)$, $\hat{y}(x) = \argmax\limits_{y \in \mathcal{Y}}\hat{F}(x, y)$.

We require the standard strong convexity and strong concavity under (arbitrary) norms of the objective function $F(x, y, \xi)$ by arguments $x$ and $y$ respectively.

\begin{assumption}\label{scsc}
    For almost every $\xi \sim P$, $F(\cdot, y, \xi)$ is $\sigma_x$-strongly convex under norm $\|\cdot\|_x$. $F(x, \cdot, \xi)$ is $\sigma_y$-strongly concave under norm $\|\cdot\|_y$, i.e.,
    \begin{align}
    &\forall x, x' \in \mathcal{X}, y, y' \in \mathcal{Y}\\\nonumber
    &F(\alpha x + (1 - \alpha)x', y, \xi) \leq  \\\nonumber
    &\alpha F(x, y, \xi) + (1 - \alpha) F(x', y, \xi) - \frac{\sigma_x\alpha(1 - \alpha)}{2}\|x - x'\|_x^2,\\\nonumber
    &F(x, \alpha y + (1 - \alpha)y', \xi) \geq \\\nonumber
    &\alpha F(x, y, \xi) + (1 - \alpha) F(x, y', \xi) + \frac{\sigma_x\alpha(1 - \alpha)}{2}\|y - y'\|_y^2.
\end{align}
\end{assumption}

See \cite{kakade2012regularization}, \cite{gonen2017average}, \cite{shalev2007online} references for importance non-Euclidean norms for regularization and preconditioning.

We further assume the function is Lipschitz continuous.
\begin{assumption}\label{lip}
     For almost every $\xi \sim P$, $F(\cdot, y, \xi)$ is $L_x$-Lipschitz continuous under norm $\|\cdot\|_x$. $F(x, \cdot, \xi)$ is $L_y$-Lipschitz continuous under norm $\|\cdot\|_y$, i.e.,
     \begin{align}
        &\forall x, x' \in \mathcal{X}, y, y' \in \mathcal{Y},\\\nonumber
        &|F(x, y, \xi) - F(x', y, \xi)| \leq L_x\|x - x'\|_x,\\\nonumber
        &|F(x, y, \xi) - F(x, y', \xi)| \leq L_x\|y - y'\|_y.
     \end{align}
\end{assumption}
Also we assume the gradient is Lipschitz continuous.
\begin{assumption}\label{lipgrad}
     For almost every $\xi \sim P$, $\nabla_x F(x, \cdot, \xi)$ and $\nabla_y F(x, \cdot, \xi)$ are $L_{x, y}$-Lipschitz continuous under norms $\|\cdot\|_x$, $\|\cdot\|_y$ respectively, i.e.,
     \begin{align}
        &\forall x, x' \in \mathcal{X}, y, y' \in \mathcal{Y}\\\nonumber
        &\|\nabla_x F(x, y, \xi) - \nabla_x F(x, y', \xi)\|_{*y} \leq L_{x, y}\|y - y'\|_y\\\nonumber
        &\|\nabla_y F(x, y, \xi) - \nabla_y F(x', y, \xi)\|_{*x} \leq L_{x, y}\|x - x'\|_x,
     \end{align}
     where $\|\cdot\|_{*x}, \|\cdot\|_{*y}$ stand for the dual norms of $\|\cdot\|_x$ and $\|\cdot\|_y$.
\end{assumption}

For our result we will need to bound the Lipschitz constant of the gradient from above with strong convexity constants.

\begin{assumption}\label{lipgradbound}
    $L_{x,y} \leq \min(\sigma_x, \sigma_y)$.
\end{assumption}

The assumption \ref{lipgradbound} can be achieved by adding strongly convex regularizer functions of $x$ and $y$, e.g., $\|x\|_x^2$ and $\|y\|_y^2$ with suitable constants.

We are ready to present our main result.

\begin{theorem}\label{maintheorem}
    Assume $\mathcal{X}, \mathcal{Y}$ to be convex, compact subsets of $\mathcal{R}^d$, under assumptions \ref{scsc}, \ref{lip}, \ref{lipgrad}, \ref{lipgradbound} for any $\delta \in (0, 1)$, it holds that, with probability at least $1 - \delta$ the ESP problem \eqref{esp} solution satisfies
    \begin{equation}
        F(\hat{x}, y^*(\hat{x})) - F(x^*(\hat{y}), \hat{y}) =  \mathcal{O}\left(\frac{d + \log(\nicefrac{1}{\delta})}{n}\right),
    \end{equation}
    where $\mathcal{O}$-constant depends only from $L_x, L_y, L_{x, y}, \sigma_x, \sigma_y$.
\end{theorem}

\section{Examples}
The minimax setup has broad applications in
mathematics and machine learning. Here we give some examples to illustrate our result.

\subsection{Matrix games}
Consider the regularized two-player stochastic matrix game problem:

\begin{equation}\label{game}
    \min_{x \in \Delta_{d}}\max_{y \in \Delta_{d}} x^T \EE_{\xi}A_{\xi}y + V_x(x) - V_y(y),
\end{equation}
where 
\begin{align}\nonumber
    \Delta_N = \{x \geq e^{-L} \mid \sum\limits_1^N x_i = 1\},
    V_x(x) = \lambda_x \sum\limits_1^{d}x_i\log(x_i),
\end{align}
and similarly defined $V_y(y)$. Functions $V_x(x)$ and $V_y(y)$ are Lipschitz continuous with constants $\lambda_x(L + 1), \lambda_y(L + 1)$ and strongly convex with constants $\lambda_x, \lambda_y$ accordingly. Assume the matrix elements are bounded: $\max_{i, j} |A_{\xi}(i, j)| \leq 1$

The solution $x^*, y^*$ of the problem \eqref{game} is the Nash equilibrium point, meaning that if the player 1
plays strategy $x^*$, it does not gain much benefit for player 2 if he switches to any other strategy.

Let us verify the assumptions: \eqref{scsc} holds with $\sigma_x = \lambda_x$, $\sigma_y = \lambda_y$, \eqref{lip} holds with $L_x = \lambda_x(L + 1) + 1$, $L_y = \lambda_y(L + 1) + 1$, \eqref{lipgrad} holds with $L_{x, y} = 1$, \eqref{lipgradbound} holds if $\lambda_x > 1$ and $\lambda_y > 1$.

By directly applying theorem \ref{maintheorem} we conclude, what solution to the empirical counterpart of problem \eqref{game} $\hat{x}, \hat{y}$ satisfies: 
\begin{align}
\nonumber
    &\hat{x}^TA_{\xi}\hat{y} \leq {(x^*)}^TA_{\xi}y^* + \mathcal{O}\left(\frac{d + \log(\nicefrac{1}{\delta})}{n}\right) + \\\nonumber
    &\mathcal{O}\left(\max_{x \in \Delta_{d}}(|V_x(x)|) + \max_{y \in \Delta_{d}} (|V_y(y)|)\right) \leq\\\nonumber 
    &{(x^*)}^TA_{\xi}y^* +
    \mathcal{O}\left(\log(d)\right). 
\end{align}
and 
\begin{align}
\nonumber
    &{(x^*)}^TA_{\xi}y^* - \mathcal{O}\left(\log(d)\right) - \mathcal{O}\left(\frac{(L + 1)(d + \log(\nicefrac{1}{\delta}))}{n}\right) \leq \hat{x}^TA_{\xi}\hat{y},
\end{align}
with probability $\geq 1 - \delta$.
Compare the last result with \cite{zhang2021generalization}, \cite{nemirovski2009robust}.

\begin{remark}
    Note that in contrast to \cite{zhang2021generalization} we solve the problem on the truncated simplex $\Delta_N$ to preserve Lipshitz continuity of the regularizer
\end{remark}

\subsection{AUC maximization}

Consider the linear version of the AUC optimization problem from \cite{ying2016stochastic}:

\begin{equation}\label{auc}
    \min_{w, a, b}\max_{\alpha \in \mathcal{R}} \left[F_{x, y}(w, a, b, \alpha, x, y) + \beta\|w\|_2^2\right],
\end{equation}
where $F(w, a, b, \alpha, x, y) = (1 - p)(w^Tx - a)^2 I[y=1] + p(w^Tx - b)^2I[y=-1] + 2(1 + \alpha)w^Tx(pI[y=-1] - (1 - p)I[y=1]) - p(1 - p)\alpha^2$, $x \in B_r \subset \mathcal{R}^d, y \in \{+1, -1\}$ are drawn from sample data distribution, $w \in B_r \subset \mathcal{R}^d, a, b \in \mathcal{R}$ --- parameters of the linear model, $\alpha \in \mathcal{R}$, $B_r$ is the ball with radius $r$ and center in point $0$.

The function in the problem \eqref{auc} is Lipschitz with $L_x = L_y = \mathcal{O}(\beta r)$, strongly convex with $\sigma_x = \sigma_y = \mathcal{O}(\beta)$, and it holds \eqref{lipgrad} with $L_{x, y} = \mathcal{O}(r)$. Applying theorem \ref{maintheorem} we get $\hat{w}, \hat{a}, \hat{b}$ as the model parameters with AUC score, which is different from different from the 
optimum by $\mathcal{O}(\frac{\beta r(d + \log\nicefrac{1}{\delta})}{\beta n})$ with probability $\geq 1 - \delta$.

Compare the last result with \cite{ying2016stochastic}.

\section{Proof}
%мы развиваем
Our proof relies on the standard arguments on the Laplace transform of shifted or offset empirical processes. We develop a similar approach was proposed in \cite{zhivotovskiy2018localization}, \cite{kanade2022exponential}, \cite{puchkin2023exploring} for minimization problems. We naturally generalize it for SSP problems.

Major differences from previous works:
\begin{enumerate}
    \item In contrast to \cite{puchkin2023exploring}, where local norms are explored in the minimization problem, we consider the saddle-point problem with non-separable variables and propose a new way of localization (Lemma \ref{lemma_42}, appendix \ref{appendix2}).
    \item We prove the possibility of localization relatively to saddle-point $x^*, y^*$ and get strong risk bound, compared to localization relatively to $x^*(\hat{y}), y^*(\hat{x})$, which can be obtained by a natural generalization of \cite{puchkin2023exploring} to saddle-point problems and leads to weak risk bound.
\end{enumerate}

\begin{lemma}\label{lemma_41}
    Given $\Phi: \R \rightarrow \R$ convex monotonously increasing function. Under assumptions \ref{scsc}, \ref{lip}, \ref{lipgrad} holds:
    \begin{align}\label{lemma_41_bound}
         &\EE_{S_n}\left[\Phi\left(F(\hat{x}, y^*(\hat{x})) - F(x^*(\hat{y}), \hat{y})\right)\right] \leq\\\nonumber &\EE_{S_n}\EE_{\e}\left[\Phi\left(4 \sup\limits_{y \in \mathcal{Y}, x \in \mathcal{X}}\left[\PP_n\e(F(x, y^*(x),\xi) - F(x^*(y),y,\xi)) -\right.\right.\right.\\\nonumber
         &\left.\left.\left.\frac{\sigma_y}{8}\|y - y^*(x)\|_y^2 - \frac{\sigma_x}{8}\|x - x^*(y)\|_x^2 ))\vphantom{\sup\limits_{y \in \mathcal{Y}, x \in \mathcal{X}}}\right]\right)\right],
    \end{align}
    where $\e_1, ..., \e_n$ are i.i.d. Rademacher random variables.
\end{lemma}

The proof is based on the one presented in Appendix A.3. \cite{puchkin2023exploring} and is given in \ref{appendix1}.

Let us focus on the shifted multiplier process of the form $\PP_n\e(F(x, y^*(x),\xi) - F(x^*(y),y,\xi)) - \frac{\sigma_y}{8}\|y - y^*(x)\|_y^2 - \frac{\sigma_x}{8}\|x - x^*(y)\|_x^2 ))$ and consider exponential moments of its supremum.

\begin{lemma}\label{lemma_42}
    Set $\lambda = \frac{\max(\sigma_x, \sigma_y)\left(1 - \frac{L_{x,y}}{\min(\sigma_x, \sigma_y)}\right)^2n}{32\sqrt{2}e\left(2\max(L_x, L_y)\sqrt{2}\left(1 + \frac{L_{x,y}}{\min(\sigma_x, \sigma_y)}\right)\right)^2}$. Under assumptions \ref{scsc}, \ref{lip}, \ref{lipgrad}, \ref{lipgradbound}. For any sample $S_n$ holds.
    \begin{align}\label{lemma_42_bound}
        &\EE_{\e}\left[\text{exp}\left\{\lambda\sup\limits_{y \in \mathcal{Y}, x \in \mathcal{X}}\left(\PP_n\e(F(x, y^*(x),\xi) - F(x^*(y),y,\xi)) -\right.\right.\right.\\\nonumber
        &\left.\left.\left.\frac{\sigma_y}{8}\|y - y^*(x)\|_y^2 - \frac{\sigma_x}{8}\|x - x^*(y)\|_x^2))\vphantom{\sup\limits_{y \in \mathcal{Y}, x \in \mathcal{X}}}\right)\right\}\right] \leq\\\nonumber
        &e + e^{3d} + 12e^{2048(1 + e)^2\nicefrac{d}{e}}.
    \end{align}
\end{lemma}

Lemmas \ref{lemma_41} and \ref{lemma_42} imply
\begin{equation}
    \log\EE\text{exp}\{\lambda\left(F(\hat{x}, y^*(\hat{x})) - F(x^*(\hat{y}), \hat{y})\right)\} = \mathcal{O}(d).
\end{equation}
Applying Markov inequality we get
\begin{equation}
    F(\hat{x}, y^*(\hat{x})) - F(x^*(\hat{y}), \hat{y}) = \mathcal{O}\left(\frac{d + \log(\nicefrac{1}{\delta})}{n}\right).
\end{equation}

\section{Conclusion}

In this work we have proved the tight high-probability upper bound on the excess risk of ERM for strongly convex-strongly concave losses with Lipschitz continuous gradient with appropriate regularization. We have generalized the local-norms analysis from \cite{puchkin2023exploring} of the minimax setup, but our result doesn't hold for exp-concave-exp-convex losses. The question on existing $\mathcal{O}\left(\frac{d + \log(\nicefrac{1}{\delta})}{n}\right)$ high-probability bound for ESP solution for exp-concave setup remains open, as well as the question about relaxing the Assumption \eqref{lipgradbound}. In addition, it is interesting to prove similar generalization bound for online setup (gradient-type algorithms).

\section{Acknowledgements}

The author would like to thank Alexander Gasnikov and Nikita Puchkin for their valuable advice and assistance in compiling the paper.

\begin{comment}
\section*{Declarations}

\subsection*{Funding}
No funding was received to assist with the preparation of this manuscript.
\subsection*{Competing interests}
The authors have no relevant financial or non-financial interests to disclose.
\subsection*{Ethics approval}
Not applicable.
\subsection*{Consent to participate}
Not applicable.
\subsection*{Consent for publication}
All authors whose names appear on the submission made substantial contributions to the conception or design of the work, approved the version to be published and agree to be accountable for all aspects of the work in ensuring that questions related to the accuracy or integrity of any part of the work are appropriately investigated and resolved.
\subsection*{Availability of data and materials}
Not applicable.
\subsection*{Code availability}
Not applicable.
\subsection*{Authors' contributions}
Vitali Pirau has written the whole text.
\end{comment}
\begin{appendices}
\section{Proof of Lemma \ref{lemma_41}}\label{appendix1}

\begin{proof}
    From \eqref{scsc} get:
    \begin{align}
    \nonumber
        &\PP\left(F(\hat{x}, \hat{y},\xi) - F(\hat{x}, y^*(\hat{x}),\xi)\right) \leq -\frac{\sigma_y}{4}\|\hat{y} - y^*(\hat{x})\|_{y}^2\\\nonumber
        &\PP_n\left(F(\hat{x}, y^*(\hat{x}),\xi) - F(\hat{x}, \hat{y},\xi)\right) \leq -\frac{\sigma_y}{4}\|\hat{y} - y^*(\hat{x})\|_{x}^2.
    \end{align}
    Note that
    \begin{align}\label{ssp_rewrite}
         F(\hat{x}, y^*(\hat{x})) - F(\hat{x}, \hat{y}) =  &2(\PP - \PP_n)\left(F(\hat{x}, y^*(\hat{x}),\xi) - F(\hat{x}, \hat{y},\xi)\right) +\\\nonumber
         &\PP\left(F(\hat{x}, \hat{y},\xi) - F(\hat{x}, y^*(\hat{x}),\xi)\right) +2\PP_n\left(F(\hat{x}, y^*(\hat{x}),\xi) - F(\hat{x}, \hat{y},\xi)\right) \le\\\nonumber
         &2(\PP - \PP_n)\left(F(\hat{x}, y^*(\hat{x}),\xi) - F(\hat{x}, \hat{y},\xi)\right) -\frac{\sigma_y}{4}\|\hat{y} - y^*(\hat{x})\|_y^2 - \frac{\sigma_y}{2}\|\hat{y} - y^*(\hat{x})\|_{y}^2.
    \end{align}
    Similarly,
    \begin{align}\label{esp_rewrite}
         F(\hat{x}, \hat{y}) - F(x^*(\hat{y}), \hat{y})=&2(\PP - \PP_n)\left(F(\hat{x}, \hat{y},\xi) - F(x^*(\hat{y}), \hat{y},\xi)\right) +\\ \nonumber
         &\PP\left(F(x^*(\hat{y}), \hat{y},\xi) - F(\hat{x}, \hat{y},\xi)\right) + 2\PP_n\left(F(\hat{x}, \hat{y},\xi) - F(x^*(\hat{y}), \hat{y},\xi)\right) \le\\\nonumber
         &2(\PP - \PP_n)\left(F(\hat{x}, \hat{y},\xi) - F(x^*(\hat{y}), \hat{y},\xi)\right) - \frac{\sigma_x}{4}\|\hat{x} - x^*(\hat{y})\|_{x}^2 - \frac{\sigma_y}{2}\|\hat{y} - y^*(\hat{x})\|_{y}^2.
    \end{align}
    Adding the inequalities \eqref{ssp_rewrite} and \eqref{esp_rewrite}, get
    \begin{align}\label{monotonic_emp_bound}
        &F(\hat{x}, y^*(\hat{x})) - F(x^*(\hat{y}), \hat{y})= \\\nonumber
        &2(\PP - \PP_n)\left(F(\hat{x}, y^*(\hat{x}),\xi)- F(x^*(\hat{y}), \hat{y},\xi)\right) - \frac{3\sigma_y}{4}\|\hat{y} - y^*(\hat{x})\|_{y}^2 -  \frac{3\sigma_x}{4}\|\hat{x} - x^*(\hat{y})\|_{x}^2.
    \end{align}
    Now we will apply standard simmetrization argument with i.i.d. Rademacher random values $\{\e_i\}_{i = 1}^n$. 
    Let us bound \eqref{monotonic_emp_bound} with the expression below.
    \begin{align}
    \nonumber
        2\sup\limits_{y \in \mathcal{Y}, x \in \mathcal{X}}
        &\left(\vphantom{\frac{\sigma_x}{4}}(\PP - \PP_n)\left(F(x, y^*(x),\xi) - F(x^*(y), y, \xi)\right) -\right.\\\nonumber
        &\left.\frac{\sigma_y}{4}\|y - y^*(x)\|_{y}^2 - \frac{\sigma_x}{4}\|x - x^*(y)\|_{x}^2\right).
    \end{align}
    We introduce independent copy of $S_n$, $S'_n = \{\xi_1', ..., \xi_n'\}$, $\PP_n'$ --- expectation with respect to uniform measure on $S'_n$, $\EE'_{S'_n}$ --- conditional expectation with respect $S_n'$.
    \begin{align}\label{symmetry_argument}
        &\EE_{S_n}\left[\Phi\left(2\sup_{y \in \mathcal{Y}, x \in \mathcal{X}}\left((\PP - \PP_n)\left(F(x, y^*(x),\xi) - F(x^*(y), y,\xi)\right) - \frac{\sigma_y}{4}\|y - y^*(x)\|_{y}^2 - \frac{\sigma_x}{4}\|x - x^*(y)\|_{x}^2\right)\right)\right] =\\\nonumber
        &\EE_{S_n}\left[\Phi\left(2\sup_{y \in \mathcal{Y}, x \in \mathcal{X}}\left(\EE'_{S'_n}(\PP_n' - \PP_n)\left(F(x, y^*(x),\xi) - F(x^*(y), y,\xi)\right)\right.\right. \right.  \\\nonumber
        &\left.\left.\left.- \frac{\sigma_y}{4}\|y - y^*(x)\|_{y}^2 - \frac{\sigma_x}{4}\|x - x^*(y)\|_{x}^2\right)\right)\right] \leq \\\nonumber
    &\EE_{S_n}\EE'_{S'_n}\left[\Phi\left(2\sup_{y \in \mathcal{Y}, x \in \mathcal{X}}\left((\PP_n' - \PP_n)\left(F(x, y^*(x),\xi) - F(x^*(y), y,\xi)\right) \right.\right. \right.  \\\nonumber
        &\left.\left.\left.- \frac{\sigma_y}{4}\|y - y^*(x)\|_{y}^2 - \frac{\sigma_x}{4}\|x - x^*(y)\|_{x}^2\right)\right)\right] = \\\nonumber
        &\EE_{S_n}\EE'_{S'_n}\EE_{\e}\left[\Phi\left(2\sup_{y \in \mathcal{Y}, x \in \mathcal{X}}\left((\PP_n' - \PP_n)\e\left(F(x, y^*(x),\xi) - F(x^*(y), y,\xi)\right)  \right.\right. \right.  \\\nonumber
        &\left.\left.\left. - \frac{\sigma_y}{4}\|y - y^*(x)\|_{y}^2- \frac{\sigma_x}{4}\|x - x^*(y)\|_{x}^2\right)\right)\right] \le\\\nonumber
        &\EE_{S_n}\EE_{\e}\left[\Phi\left(4\sup_{y \in \mathcal{Y}, x \in \mathcal{X}}\left(\PP_n\e\left(F(x, y^*(x),\xi) - F(x^*(y), y,\xi)\right) \right.\right. \right.  \\\nonumber
        &\left.\left.\left.- \frac{\sigma_y}{8}\|y - y^*(x)\|_{y}^2 - \frac{\sigma_x}{8}\|x - x^*(y)\|_{x}^2\right)\right)\right].
    \end{align}
    Substitute \eqref{symmetry_argument} in \eqref{monotonic_emp_bound}
    \begin{align}
         \nonumber
         &\EE_{S_n}\left[\Phi\left(F(\hat{x}, y^*(\hat{x})) - F(x^*(\hat{y}), \hat{y})\right)\right] \leq\\\nonumber &\EE_{S_n}\EE_{\e}\left[\Phi\left(4 \sup\limits_{y \in \mathcal{Y}, x \in \mathcal{X}}\left[\PP_n\e(F(x, y^*(x),\xi) - F(x^*(y),y,\xi)) - \frac{\sigma_y}{8}\|y - y^*(x)\|_y^2 - \frac{\sigma_x}{8}\|x - x^*(y)\|_x^2 ))\right]\right)\right].
    \end{align}
\end{proof}

\section{Proof of Lemma \ref{lemma_42}}\label{appendix2}

\begin{proof}
    Let us consider arbitrary value of the sample $S_n$ fixed.
    Denote $\mathcal{XY}[a, b] = \{(x, y) \in \mathcal{X}\times\mathcal{Y} \mid a \le \|y - y^*\|_y + \|x - x^*\|_x \le b\}$.
    \begin{align}\label{localbound}
        &\EE_{\e}\left[\text{exp}\left\{\lambda\sup\limits_{y \in \mathcal{Y}, x \in \mathcal{X}}\left(\PP_n\e(F(x, y^*(x),\xi) - F(x^*(y),y,\xi)) - \frac{\sigma_y}{8}\|y - y^*(x)\|_y^2 - \frac{\sigma_x}{8}\|x - x^*(y)\|_x^2))\right)\right\}\right] \le \\\nonumber
        &\EE_{\e}\left[\text{exp}\left\{\lambda\sup\limits_{(x, y) \in \mathcal{XY}[0, r]}\left(\PP_n\e(F(x, y^*(x),\xi) - F(x^*(y),y,\xi)) - \frac{\sigma_y}{8}\|y - y^*(x)\|_y^2 - \frac{\sigma_x}{8}\|x - x^*(y)\|_x^2))\right)\right\}\right] +\\\nonumber
        &\sum\limits_{k = 0}^\infty\EE_{\e}\left[\text{exp}\left\{\lambda\sup\limits_{(x, y) \in \mathcal{XY}[2^{k}r, 2^{k+1}r]}\left(\PP_n\e(F(x, y^*(x),\xi) - F(x^*(y),y,\xi))\right.\right. \right.  \\\nonumber
        &\left.\left.\left. - \frac{\sigma_y}{8}\|y - y^*(x)\|_y^2 - \frac{\sigma_x}{8}\|x - x^*(y)\|_x^2))\right)\right\}\right].
    \end{align}
    Now we will bound the expectation of the exponent of the local supremum from above.
    For this purpose, we will bound the negative terms in the supremum from below.
    From \ref{lipgrad} and  \cite[Lemma 4, Section A2]{zhang2021generalization} get:
    \begin{align}\label{gradnormchange}
        \|x^*(y) - x^*\|_{x} \leq \frac{L_{x,y}}{\sigma_x}\|y - y^*\|_{y}\\\nonumber
        \|y^*(x) - y^*\|_{y} \leq \frac{L_{x,y}}{\sigma_y}\|x - x^*\|_{x}
    \end{align}
    From triangle inequality:
    \begin{align}
    \nonumber
        \|x - x^*(y)\|_{x} \geq \|x - x^*\|_{x} - \|x^* - x^*(y)\|_{x} \geq \|x - x^*\|_{x} - \frac{L_{x,y}}{\sigma_x}\|y^* - y\|_{y}\\\nonumber
        \|y - y^*(x)\|_{y} \geq \|y - y^*\|_{y} - \frac{L_{x,y}}{\sigma_y}\|x^* - x\|_{x}
    \end{align}
    From \eqref{gradnormchange}:
    \begin{align}
    \nonumber
         \left(1 - \frac{L_{x,y}}{\sigma_y}\right)\|x - x^*\|_{x} + \left(1 - \frac{L_{x,y}}{\sigma_x}\right)\|y - y^*\|_{y} \leq \|x - x^*(y)\|_{x} + \|y - y^*(x)\|_{y},
    \end{align} and
    \begin{align}\label{regularizer_change}
         \frac{C^2}{2}\left(\|x - x^*\|_{x} + \|y - y^*\|_{y}\right)^2 \leq \|x - x^*(y)\|^2_{x} + \|y - y^*(x)\|^2_{y}.
    \end{align}
    where $C = \left(1 - \frac{L_{x,y}}{\min(\sigma_x, \sigma_y)}\right)$.
    We generalize \cite[Lemma 5.1]{puchkin2023exploring} on SSP problems and find an upper bound on the exponential moment of the supremum of a localized set.
    %\esssup
    \begin{lemma}\label{lemma_51}
        Denote:
        \begin{align}
        \nonumber
            &B = \sup\limits_{(x, y) \in \mathcal{XY}\left[0, r\right]}\esssup\limits_{\xi \in P}|F(x, y^*(x),\xi) - F(x^*(y), y,\xi)|\\\nonumber
            &L = \max(L_x, L_y)\\\nonumber
            &\Tilde{L} = 2L\Tilde{C}\\\nonumber
            &\Tilde{C} = \sqrt{2}\left(1 + \frac{L_{x,y}}{\min(\sigma_x, \sigma_y)}\right).
        \end{align}
        Under assumptions \ref{scsc}, \ref{lip}, \ref{lipgrad}, \ref{lipgradbound} for all $\lambda \ge 0$ holds
        \begin{align}
        \nonumber
            &\log\EE_{\e}\left[\text{exp}\left\{\lambda\sup\limits_{(x, y) \in \mathcal{XY}\left[0, r\right]}\left(\PP_n\e(F(x, y^*(x),\xi)-F(x^*(y), y,\xi)\right)\right\}\right] \leq\\ &64\sqrt{2}\lambda \Tilde{L} r \sqrt{\frac{d}{n}} + \frac{B^2\lambda^2e^{\nicefrac{B\lambda}{n}}\sqrt{2}}{2n}\left(\frac{128\Tilde{L}r}{B}\sqrt{\frac{d}{n}} + \frac{\Tilde{L}^2 r^2}{B^2}\right).
        \end{align}
        If, in addition $\lambda \le \frac{n}{B}$, then
        \begin{align}
        \nonumber
            &\log\EE_{\e}\left[\text{exp}\left\{\lambda\sup\limits_{(x, y) \in \mathcal{XY}\left[0, r\right]}\left(\PP_n\e(F(x, y^*(x),\xi)-F(x^*(y), y,\xi)\right)\right\}\right] \leq 64\sqrt{2}(1+e)\lambda \Tilde{L} r \sqrt{\frac{d}{n}} + \frac{e\lambda^2 \Tilde{L}^2 r^2\sqrt{2}}{2n}.
        \end{align}
    \end{lemma}
    The Lemma \ref{lemma_51} is proved by explicit construction of an epsilon net.

    From Assumption \ref{lip}:
    \begin{align}
    \nonumber
        \forall x \in \mathcal{X}, y_1, y_2 \in \mathcal{Y} \; \PP_n\left(F(x, y_1, \xi) - F(x, y_2, \xi)\right)^2 \leq L_y^2\|y_1 - y_2\|_H^2\\\nonumber
        \forall y \in \mathcal{Y}, x_1, x_2 \in \mathcal{X} \; \PP_n\left(F(x_1, y, \xi) - F(x_2, y, \xi)\right)^2 \leq L_x^2\|x_1 - x_2\|_H^2.
    \end{align}
    Denote $\mathcal{F}[0, r] = \{F(x, y^*(x), \xi) - F(x^*(y), y, \xi) \mid (x, y) \in \mathcal{XY}\left[0, r\right]\}$. Then a pair of $\nicefrac{\gamma}{2L\Tilde{C}}$-nets on balls in $\mathcal{X}, \mathcal{Y}$ with centers $x^*, y^*$ accordingly and radii $r$  generates a $\gamma$-net on $\mathcal{F}[0, r]$ in empirical $L_2$-norm. Let $x_{\nicefrac{\gamma}{L}}$ be the closest to $x$ element in epsilon net on $\mathcal{X}$-ball, $y_{\nicefrac{\gamma}{L}}$ the closest to $y$ in epsilon net on $\mathcal{Y}$-ball. We will show that $F(x_{\nicefrac{\gamma}{L}}, y^*(x_{\nicefrac{\gamma}{L}}), \xi) - F(x^*(y_{\nicefrac{\gamma}{L}}), y_{\nicefrac{\gamma}{L}}, \xi) \in B_{\gamma}\left(F(x, y^*(x), \xi) - F(x^*(y), y, \xi)\right)$ in $\mathcal{F}[0, r]$. The distance in the empirical $L_2$-norm is:
    \begin{align}
        %&\left|\sqrt{\PP_n\left(F(x, y^*(x), \xi) - F(x^*(y), y, \xi)\right)^2} - \sqrt{\PP_n\left(F(x_{\nicefrac{\gamma}{L}}, y^*(x_{\nicefrac{\gamma}{L}}), \xi) - F(x^*(y_{\nicefrac{\gamma}{L}}), y_{\nicefrac{\gamma}{L}}, \xi)\right)^2}\right| \leq \\
        \nonumber
        \sqrt{\PP_n\left(F(x, y^*(x), \xi) - F(x^*(y), y, \xi) - F(x_{\nicefrac{\gamma}{L}}, y^*(x_{\nicefrac{\gamma}{L}}), \xi) + F(x^*(y_{\nicefrac{\gamma}{L}}), y_{\nicefrac{\gamma}{L}}, \xi)\right)^2}.
    \end{align}
    Under the last root, we add and subtract $F(x, y^*(x_{\nicefrac{\gamma}{L}}), \xi) + F(x^*(y_{\nicefrac{\gamma}{L}}), y, \xi)$, decompose the square of the sum into the sum of squares and apply \eqref{lip} four times.

    \begin{align}
    \nonumber
        &\sqrt{\PP_n\left(F(x, y^*(x), \xi) - F(x^*(y), y, \xi) - F(x_{\nicefrac{\gamma}{L}}, y^*(x_{\nicefrac{\gamma}{L}}), \xi) + F(x^*(y_{\nicefrac{\gamma}{L}}), y_{\nicefrac{\gamma}{L}}, \xi)\right)^2} \leq \\\nonumber
        &L\sqrt{2\left(\|y^*(x) - y^*(x_{\nicefrac{\gamma}{L}})\|_y^2 + \|x - x_{\nicefrac{\gamma}{L}}\|_x^2 + \|x^*(y) - x^*(y_{\nicefrac{\gamma}{L}})\|_x^2 + \|y - y_{\nicefrac{\gamma}{L}}\|_y^2\right)} \leq \\\nonumber
        &L\sqrt{2}\left(1 + \frac{L_{x,y}}{\min(\sigma_x, \sigma_y)}\right)\left(\|x - x_{\nicefrac{\gamma}{L}}\|_H + \|y - y_{\nicefrac{\gamma}{L}}\|_H\right) \leq \gamma.
    \end{align}
    Denote $\Tilde{L} = 2L\Tilde{C}$. Note $\mathcal{N}(\mathcal{XY}[0, r], \|\cdot\|_H, u) \leq \mathcal{N}\left(\mathcal{X}\left[0, r\right], \|\cdot\|_x, u\right)\mathcal{N}\left(\mathcal{Y}\left[0, r\right], \|\cdot\|_y, u\right)$, so we can use \cite[Lemma 5.1]{puchkin2023exploring}, with $\Tilde{L}$ as $L$ and with additional $\sqrt{2}$ factor.
    
    Continue:
    \begin{align}\label{zeroth}
        &\EE_{\e}\left[\text{exp}\left\{\lambda\sup\limits_{(x, y) \in \mathcal{XY}[0, r]}\left(\PP_n\e(F(x, y^*(x),\xi) - F(x^*(y),y,\xi)) - \frac{\sigma_y}{8}\|y - y^*(x)\|_H^2 - \frac{\sigma_x}{8}\|x - x^*(y)\|_H^2))\right)\right\}\right] \leq\\\nonumber
        &\text{exp}\left\{64\sqrt{2}(1+e)\lambda \Tilde{L} r \sqrt{\frac{d}{n}} + \frac{e\sqrt{2}\lambda^2 \Tilde{L}^2 r^2}{2n}\right\}.
    \end{align}
    Fix $k \geq 1$. Using \eqref{regularizer_change}:
    \begin{align}
    \nonumber
        &\EE_{\e}\left[\text{exp}\left\{\lambda\sup\limits_{(x, y) \in \mathcal{XY}[2^{k}r, 2^{k+1}r]}\left(\PP_n\e(F(x, y^*(x),\xi) - F(x^*(y),y,\xi)) \right.\right. \right.  \\\nonumber
        &\left.\left.\left.- \frac{1}{8}\left(\sigma_y\|y - y^*(x)\|_y^2 + \sigma_x\|x - x^*(y)\|_x^2))\right)\right)\right\}\right] \leq \\\nonumber
        &\EE_{\e}\left[\text{exp}\left\{\lambda\sup\limits_{(x, y) \in \mathcal{XY}[2^{k}r, 2^{k+1}r]}\left(\PP_n\e(F(x, y^*(x),\xi) - F(x^*(y),y,\xi)) \right.\right. \right.  \\\nonumber
        &\left.\left.\left.- \frac{\max(\sigma_x, \sigma_y)}{8}\left(\|y - y^*(x)\|_y^2 + \|x - x^*(y)\|_x^2))\right)\right)\right\}\right] \leq \\\nonumber
        &\EE_{\e}\left[\text{exp}\left\{\lambda\sup\limits_{(x, y) \in \mathcal{XY}[2^{k}r, 2^{k+1}r]}\left(\PP_n\e(F(x, y^*(x),\xi) - F(x^*(y),y,\xi)) \right.\right. \right.  \\\nonumber
        &\left.\left.\left.- \frac{C^2\max(\sigma_x, \sigma_y)}{16}\left(\|y - y^*\|_y^2 +\|x - x^*\|_x^2))\right)\right)\right\}\right].
    \end{align}
    From definition of $\mathcal{XY}[2^kr, 2^{k+1}r]$ get $\|y - y\|_y + \|x - x^*\|_x \geq 2^kr$. Then we get:
    \begin{align}\label{kth}
        &\EE_{\e}\left[\text{exp}\left\{\lambda\sup\limits_{(x, y) \in \mathcal{XY}[2^{k}r, 2^{k+1}r]}\left(\PP_n\e(F(x, y^*(x),\xi) - F(x^*(y),y,\xi)) \right.\right. \right.  \\\nonumber
        &\left.\left.\left.- \frac{1}{8}\left(\sigma_y\|y - y^*(x)\|_y^2 + \sigma_x\|x - x^*(y)\|_x^2))\right)\right)\right\}\right] \leq \\\nonumber
        &\EE_{\e}\left[\text{exp}\left\{\lambda\sup\limits_{(x, y) \in \mathcal{XY}[2^{k}r, 2^{k+1}r]}\left(\PP_n\e(F(x, y^*(x),\xi) - F(x^*(y),y,\xi))\right) - C^2\max(\sigma_x, \sigma_y)4^{k-2}\lambda r^2\right\}\right] \leq \\\nonumber
        &\EE_{\e}\left[\text{exp}\left\{\lambda\sup\limits_{(x, y) \in \mathcal{XY}[0, 2^{k+1}r]}\left(\PP_n\e(F(x, y^*(x),\xi) - F(x^*(y),y,\xi))\right) - C^2\max(\sigma_x, \sigma_y)4^{k-2}\lambda r^2\right\}\right] \le \\\nonumber
        &\text{exp}\left\{64(1+e)\lambda \Tilde{L} 2^{k+1}r \sqrt{\frac{d}{n}} + \frac{e\lambda^2 \Tilde{L}^2 4^{k+1}r^2}{2n} - C^2\max(\sigma_x, \sigma_y)4^{k-2}\lambda \nicefrac{r^2}{2}\right\}.
    \end{align}
    The last part of the proof follows proof of \cite[Lemma 4.2]{puchkin2023exploring} with a little changes.
    Recall $\lambda = \frac{C^2\max(\sigma_x, \sigma_y)n}{32\sqrt{2}e\Tilde{L}^2}$, so
    \begin{equation}
    \nonumber
        \frac{e\sqrt{2}\lambda^2 \Tilde{L}^2 4^{k+1}r^2}{2n} = \frac{2\sqrt{2}e\lambda^2 \Tilde{L}^2 4^{k}r^2}{n} = \frac{4^k\lambda r^2C^2\max(\sigma_x, \sigma_y)}{16}
    \end{equation}
    Combining \eqref{localbound}, \eqref{zeroth} and \eqref{kth} we get:
    \begin{align}\label{seriesbound}
         &\EE_{\e}\left[\text{exp}\left\{\lambda\sup\limits_{y \in \mathcal{Y}, x \in \mathcal{X}}\left(\PP_n\e(F(x, y^*(x),\xi) - F(x^*(y),y,\xi)) - \frac{\sigma_y}{8}\|y - y^*(x)\|_y^2 - \frac{\sigma_x}{8}\|x - x^*(y)\|_x^2))\right)\right\}\right] \le \\\nonumber
         &\text{exp}\left\{64\sqrt{2}(1+e)\lambda \Tilde{L} r \sqrt{\frac{d}{n}} + \frac{e\sqrt{2}\lambda^2 \Tilde{L}^2 r^2}{2n}\right\} + \\\nonumber
         &\sum\limits_{k=0}^{+\infty}\text{exp}\left\{64\sqrt{2}(1+e)\lambda \Tilde{L} 2^{k+1}r \sqrt{\frac{d}{n}} + \frac{e\sqrt{2}\lambda^2 \Tilde{L}^2 4^{k+1}r^2}{2n} - C^2\max(\sigma_x, \sigma_y)4^{k-2}\lambda r^2\right\}.   
    \end{align}
    We choose $r = \frac{\Tilde{L}}{C^2\max(\sigma_x, \sigma_y)}\sqrt{\frac{d}{n}}$ and obtain:
    \begin{align}
    \nonumber
        &64\sqrt{2}(1+e)\lambda \Tilde{L} r \sqrt{\frac{d}{n}} + \frac{e\sqrt{2}\lambda^2 \Tilde{L}^2 r^2}{2n} = \\\nonumber
        &\frac{64\sqrt{2}(1+e)\lambda d\Tilde{L}^2}{nC^2\max(\sigma_x, \sigma_y)} + \frac{\lambda\Tilde{L}^2d}{64nC^2\max(\sigma_x, \sigma_y)} \leq 3d
    \end{align}
    for the first summand element of the \eqref{seriesbound}.
    For the $k$-th series element we get:
    \begin{align}
    \nonumber
        128\sqrt{2}(1+e)\lambda \Tilde{L} 2^{k}r \sqrt{\frac{d}{n}} - C^2\max(\sigma_x, \sigma_y)4^{k-2}\lambda r^2 \leq \frac{\lambda\Tilde{L}^2d}{C^2\max(\sigma_x, \sigma_y)n}\left(128\sqrt{2}(1+e)2^k - 4^{k-2}\right).
    \end{align}
    The proof proceeds analogously to one in \cite{puchkin2023exploring} with $C^2\max(\sigma_x, \sigma_y)$ instead of $\sigma$ and $\Tilde{L}$ instead of $L$.
\end{proof}

%, $r = \frac{\Tilde{L}\sqrt{\left(\frac{2}{(\sigma_x - L_{x, y})(\sigma_y - L_{x, y})}\right)}}{\max(\sigma_x, \sigma_y)}\sqrt{\frac{d}{n}}$

By substituting \eqref{lemma_42_bound} in \eqref{lemma_41_bound} with $\Phi(x) = \text{exp}(4\lambda x)$ we get the final result.

\end{appendices}


\bibliography{references}% common bib file


\begin{thebibliography}{38}
\ifx \bisbn   \undefined \def \bisbn  #1{ISBN #1}\fi
\ifx \binits  \undefined \def \binits#1{#1}\fi
\ifx \bauthor  \undefined \def \bauthor#1{#1}\fi
\ifx \batitle  \undefined \def \batitle#1{#1}\fi
\ifx \bjtitle  \undefined \def \bjtitle#1{#1}\fi
\ifx \bvolume  \undefined \def \bvolume#1{\textbf{#1}}\fi
\ifx \byear  \undefined \def \byear#1{#1}\fi
\ifx \bissue  \undefined \def \bissue#1{#1}\fi
\ifx \bfpage  \undefined \def \bfpage#1{#1}\fi
\ifx \blpage  \undefined \def \blpage #1{#1}\fi
\ifx \burl  \undefined \def \burl#1{\textsf{#1}}\fi
\ifx \doiurl  \undefined \def \doiurl#1{\url{https://doi.org/#1}}\fi
\ifx \betal  \undefined \def \betal{\textit{et al.}}\fi
\ifx \binstitute  \undefined \def \binstitute#1{#1}\fi
\ifx \binstitutionaled  \undefined \def \binstitutionaled#1{#1}\fi
\ifx \bctitle  \undefined \def \bctitle#1{#1}\fi
\ifx \beditor  \undefined \def \beditor#1{#1}\fi
\ifx \bpublisher  \undefined \def \bpublisher#1{#1}\fi
\ifx \bbtitle  \undefined \def \bbtitle#1{#1}\fi
\ifx \bedition  \undefined \def \bedition#1{#1}\fi
\ifx \bseriesno  \undefined \def \bseriesno#1{#1}\fi
\ifx \blocation  \undefined \def \blocation#1{#1}\fi
\ifx \bsertitle  \undefined \def \bsertitle#1{#1}\fi
\ifx \bsnm \undefined \def \bsnm#1{#1}\fi
\ifx \bsuffix \undefined \def \bsuffix#1{#1}\fi
\ifx \bparticle \undefined \def \bparticle#1{#1}\fi
\ifx \barticle \undefined \def \barticle#1{#1}\fi
\bibcommenthead
\ifx \bconfdate \undefined \def \bconfdate #1{#1}\fi
\ifx \botherref \undefined \def \botherref #1{#1}\fi
\ifx \url \undefined \def \url#1{\textsf{#1}}\fi
\ifx \bchapter \undefined \def \bchapter#1{#1}\fi
\ifx \bbook \undefined \def \bbook#1{#1}\fi
\ifx \bcomment \undefined \def \bcomment#1{#1}\fi
\ifx \oauthor \undefined \def \oauthor#1{#1}\fi
\ifx \citeauthoryear \undefined \def \citeauthoryear#1{#1}\fi
\ifx \endbibitem  \undefined \def \endbibitem {}\fi
\ifx \bconflocation  \undefined \def \bconflocation#1{#1}\fi
\ifx \arxivurl  \undefined \def \arxivurl#1{\textsf{#1}}\fi
\csname PreBibitemsHook\endcsname

%%% 1
\bibitem[\protect\citeauthoryear{Zhang et~al.}{2021}]{zhang2021generalization}
\begin{bchapter}
\bauthor{\bsnm{Zhang}, \binits{J.}},
\bauthor{\bsnm{Hong}, \binits{M.}},
\bauthor{\bsnm{Wang}, \binits{M.}},
\bauthor{\bsnm{Zhang}, \binits{S.}}:
\bctitle{Generalization bounds for stochastic saddle point problems}.
In: \bbtitle{International Conference on Artificial Intelligence and
  Statistics},
pp. \bfpage{568}--\blpage{576}
(\byear{2021}).
\bcomment{PMLR}
\end{bchapter}
\endbibitem

%%% 2
\bibitem[\protect\citeauthoryear{Farnia and Ozdaglar}{2021}]{farnia2021train}
\begin{bchapter}
\bauthor{\bsnm{Farnia}, \binits{F.}},
\bauthor{\bsnm{Ozdaglar}, \binits{A.}}:
\bctitle{Train simultaneously, generalize better: Stability of gradient-based
  minimax learners}.
In: \bbtitle{International Conference on Machine Learning},
pp. \bfpage{3174}--\blpage{3185}
(\byear{2021}).
\bcomment{PMLR}
\end{bchapter}
\endbibitem

%%% 3
\bibitem[\protect\citeauthoryear{Kang et~al.}{2022}]{kang2022stability}
\begin{botherref}
\oauthor{\bsnm{Kang}, \binits{Y.}},
\oauthor{\bsnm{Liu}, \binits{Y.}},
\oauthor{\bsnm{Li}, \binits{J.}},
\oauthor{\bsnm{Wang}, \binits{W.}}:
Stability and generalization of differentially private minimax problems.
arXiv preprint arXiv:2204.04858
(2022)
\end{botherref}
\endbibitem

%%% 4
\bibitem[\protect\citeauthoryear{Lei et~al.}{2021}]{lei2021stability}
\begin{bchapter}
\bauthor{\bsnm{Lei}, \binits{Y.}},
\bauthor{\bsnm{Yang}, \binits{Z.}},
\bauthor{\bsnm{Yang}, \binits{T.}},
\bauthor{\bsnm{Ying}, \binits{Y.}}:
\bctitle{Stability and generalization of stochastic gradient methods for
  minimax problems}.
In: \bbtitle{International Conference on Machine Learning},
pp. \bfpage{6175}--\blpage{6186}
(\byear{2021}).
\bcomment{PMLR}
\end{bchapter}
\endbibitem

%%% 5
\bibitem[\protect\citeauthoryear{Ji et~al.}{2021}]{ji2021understanding}
\begin{barticle}
\bauthor{\bsnm{Ji}, \binits{K.}},
\bauthor{\bsnm{Zhou}, \binits{Y.}},
\bauthor{\bsnm{Liang}, \binits{Y.}}:
\batitle{Understanding estimation and generalization error of generative
  adversarial networks}.
\bjtitle{IEEE Transactions on Information Theory}
\bvolume{67}(\bissue{5}),
\bfpage{3114}--\blpage{3129}
(\byear{2021})
\end{barticle}
\endbibitem

%%% 6
\bibitem[\protect\citeauthoryear{Goodfellow
  et~al.}{2014}]{goodfellow2014generative}
\begin{botherref}
\oauthor{\bsnm{Goodfellow}, \binits{I.}},
\oauthor{\bsnm{Pouget-Abadie}, \binits{J.}},
\oauthor{\bsnm{Mirza}, \binits{M.}},
\oauthor{\bsnm{Xu}, \binits{B.}},
\oauthor{\bsnm{Warde-Farley}, \binits{D.}},
\oauthor{\bsnm{Ozair}, \binits{S.}},
\oauthor{\bsnm{Courville}, \binits{A.}},
\oauthor{\bsnm{Bengio}, \binits{Y.}}:
Generative adversarial nets.
Advances in neural information processing systems
\textbf{27}
(2014)
\end{botherref}
\endbibitem

%%% 7
\bibitem[\protect\citeauthoryear{Du et~al.}{2017}]{du2017stochastic}
\begin{bchapter}
\bauthor{\bsnm{Du}, \binits{S.S.}},
\bauthor{\bsnm{Chen}, \binits{J.}},
\bauthor{\bsnm{Li}, \binits{L.}},
\bauthor{\bsnm{Xiao}, \binits{L.}},
\bauthor{\bsnm{Zhou}, \binits{D.}}:
\bctitle{Stochastic variance reduction methods for policy evaluation}.
In: \bbtitle{International Conference on Machine Learning},
pp. \bfpage{1049}--\blpage{1058}
(\byear{2017}).
\bcomment{PMLR}
\end{bchapter}
\endbibitem

%%% 8
\bibitem[\protect\citeauthoryear{Dai et~al.}{2018}]{dai2018sbeed}
\begin{bchapter}
\bauthor{\bsnm{Dai}, \binits{B.}},
\bauthor{\bsnm{Shaw}, \binits{A.}},
\bauthor{\bsnm{Li}, \binits{L.}},
\bauthor{\bsnm{Xiao}, \binits{L.}},
\bauthor{\bsnm{He}, \binits{N.}},
\bauthor{\bsnm{Liu}, \binits{Z.}},
\bauthor{\bsnm{Chen}, \binits{J.}},
\bauthor{\bsnm{Song}, \binits{L.}}:
\bctitle{Sbeed: Convergent reinforcement learning with nonlinear function
  approximation}.
In: \bbtitle{International Conference on Machine Learning},
pp. \bfpage{1125}--\blpage{1134}
(\byear{2018}).
\bcomment{PMLR}
\end{bchapter}
\endbibitem

%%% 9
\bibitem[\protect\citeauthoryear{Zhao et~al.}{2011}]{zhao2011online}
\begin{botherref}
\oauthor{\bsnm{Zhao}, \binits{P.}},
\oauthor{\bsnm{Hoi}, \binits{S.C.}},
\oauthor{\bsnm{Jin}, \binits{R.}},
\oauthor{\bsnm{Yang}, \binits{T.}}:
Online auc maximization
(2011)
\end{botherref}
\endbibitem

%%% 10
\bibitem[\protect\citeauthoryear{Gao et~al.}{2013}]{gao2013one}
\begin{bchapter}
\bauthor{\bsnm{Gao}, \binits{W.}},
\bauthor{\bsnm{Jin}, \binits{R.}},
\bauthor{\bsnm{Zhu}, \binits{S.}},
\bauthor{\bsnm{Zhou}, \binits{Z.-H.}}:
\bctitle{One-pass auc optimization}.
In: \bbtitle{International Conference on Machine Learning},
pp. \bfpage{906}--\blpage{914}
(\byear{2013}).
\bcomment{PMLR}
\end{bchapter}
\endbibitem

%%% 11
\bibitem[\protect\citeauthoryear{Ying et~al.}{2016}]{ying2016stochastic}
\begin{botherref}
\oauthor{\bsnm{Ying}, \binits{Y.}},
\oauthor{\bsnm{Wen}, \binits{L.}},
\oauthor{\bsnm{Lyu}, \binits{S.}}:
Stochastic online auc maximization.
Advances in neural information processing systems
\textbf{29}
(2016)
\end{botherref}
\endbibitem

%%% 12
\bibitem[\protect\citeauthoryear{Liu et~al.}{2018}]{liu2018fast}
\begin{bchapter}
\bauthor{\bsnm{Liu}, \binits{M.}},
\bauthor{\bsnm{Zhang}, \binits{X.}},
\bauthor{\bsnm{Chen}, \binits{Z.}},
\bauthor{\bsnm{Wang}, \binits{X.}},
\bauthor{\bsnm{Yang}, \binits{T.}}:
\bctitle{Fast stochastic auc maximization with $ o (1/n) $-convergence rate}.
In: \bbtitle{International Conference on Machine Learning},
pp. \bfpage{3189}--\blpage{3197}
(\byear{2018}).
\bcomment{PMLR}
\end{bchapter}
\endbibitem

%%% 13
\bibitem[\protect\citeauthoryear{Facchinei et~al.}{2014}]{facchinei2014non}
\begin{barticle}
\bauthor{\bsnm{Facchinei}, \binits{F.}},
\bauthor{\bsnm{Pang}, \binits{J.-S.}},
\bauthor{\bsnm{Scutari}, \binits{G.}}:
\batitle{Non-cooperative games with minmax objectives}.
\bjtitle{Computational Optimization and Applications}
\bvolume{59},
\bfpage{85}--\blpage{112}
(\byear{2014})
\end{barticle}
\endbibitem

%%% 14
\bibitem[\protect\citeauthoryear{Bach and Levy}{2019}]{bach2019universal}
\begin{bchapter}
\bauthor{\bsnm{Bach}, \binits{F.}},
\bauthor{\bsnm{Levy}, \binits{K.Y.}}:
\bctitle{A universal algorithm for variational inequalities adaptive to
  smoothness and noise}.
In: \bbtitle{Conference on Learning Theory},
pp. \bfpage{164}--\blpage{194}
(\byear{2019}).
\bcomment{PMLR}
\end{bchapter}
\endbibitem

%%% 15
\bibitem[\protect\citeauthoryear{Nemirovski
  et~al.}{2009}]{nemirovski2009robust}
\begin{barticle}
\bauthor{\bsnm{Nemirovski}, \binits{A.}},
\bauthor{\bsnm{Juditsky}, \binits{A.}},
\bauthor{\bsnm{Lan}, \binits{G.}},
\bauthor{\bsnm{Shapiro}, \binits{A.}}:
\batitle{Robust stochastic approximation approach to stochastic programming}.
\bjtitle{SIAM Journal on optimization}
\bvolume{19}(\bissue{4}),
\bfpage{1574}--\blpage{1609}
(\byear{2009})
\end{barticle}
\endbibitem

%%% 16
\bibitem[\protect\citeauthoryear{Zhao}{2022}]{zhao2022accelerated}
\begin{barticle}
\bauthor{\bsnm{Zhao}, \binits{R.}}:
\batitle{Accelerated stochastic algorithms for convex-concave saddle-point
  problems}.
\bjtitle{Mathematics of Operations Research}
\bvolume{47}(\bissue{2}),
\bfpage{1443}--\blpage{1473}
(\byear{2022})
\end{barticle}
\endbibitem

%%% 17
\bibitem[\protect\citeauthoryear{Shapiro et~al.}{2021}]{shapiro2021lectures}
\begin{bbook}
\bauthor{\bsnm{Shapiro}, \binits{A.}},
\bauthor{\bsnm{Dentcheva}, \binits{D.}},
\bauthor{\bsnm{Ruszczynski}, \binits{A.}}:
\bbtitle{Lectures on Stochastic Programming: Modeling and Theory}.
\bpublisher{SIAM}, \blocation{???}
(\byear{2021})
\end{bbook}
\endbibitem

%%% 18
\bibitem[\protect\citeauthoryear{Mehta}{2017}]{mehta2017fast}
\begin{bchapter}
\bauthor{\bsnm{Mehta}, \binits{N.}}:
\bctitle{Fast rates with high probability in exp-concave statistical learning}.
In: \bbtitle{Artificial Intelligence and Statistics},
pp. \bfpage{1085}--\blpage{1093}
(\byear{2017}).
\bcomment{PMLR}
\end{bchapter}
\endbibitem

%%% 19
\bibitem[\protect\citeauthoryear{Zhang et~al.}{2017}]{zhang2017empirical}
\begin{bchapter}
\bauthor{\bsnm{Zhang}, \binits{L.}},
\bauthor{\bsnm{Yang}, \binits{T.}},
\bauthor{\bsnm{Jin}, \binits{R.}}:
\bctitle{Empirical risk minimization for stochastic convex optimization: $
  o(1/n)$-and $o(1/n^2)$-type of risk bounds}.
In: \bbtitle{Conference on Learning Theory},
pp. \bfpage{1954}--\blpage{1979}
(\byear{2017}).
\bcomment{PMLR}
\end{bchapter}
\endbibitem

%%% 20
\bibitem[\protect\citeauthoryear{Wang et~al.}{2022}]{wang2022stability}
\begin{barticle}
\bauthor{\bsnm{Wang}, \binits{P.}},
\bauthor{\bsnm{Lei}, \binits{Y.}},
\bauthor{\bsnm{Ying}, \binits{Y.}},
\bauthor{\bsnm{Zhou}, \binits{D.-X.}}:
\batitle{Stability and generalization for markov chain stochastic gradient
  methods}.
\bjtitle{Advances in Neural Information Processing Systems}
\bvolume{35},
\bfpage{37735}--\blpage{37748}
(\byear{2022})
\end{barticle}
\endbibitem

%%% 21
\bibitem[\protect\citeauthoryear{Puchkin and
  Zhivotovskiy}{2023}]{puchkin2023exploring}
\begin{botherref}
\oauthor{\bsnm{Puchkin}, \binits{N.}},
\oauthor{\bsnm{Zhivotovskiy}, \binits{N.}}:
Exploring local norms in exp-concave statistical learning.
arXiv preprint arXiv:2302.10726
(2023)
\end{botherref}
\endbibitem

%%% 22
\bibitem[\protect\citeauthoryear{Natole et~al.}{2018}]{natole2018stochastic}
\begin{bchapter}
\bauthor{\bsnm{Natole}, \binits{M.}},
\bauthor{\bsnm{Ying}, \binits{Y.}},
\bauthor{\bsnm{Lyu}, \binits{S.}}:
\bctitle{Stochastic proximal algorithms for auc maximization}.
In: \bbtitle{International Conference on Machine Learning},
pp. \bfpage{3710}--\blpage{3719}
(\byear{2018}).
\bcomment{PMLR}
\end{bchapter}
\endbibitem

%%% 23
\bibitem[\protect\citeauthoryear{Yan et~al.}{2020}]{yan2020optimal}
\begin{barticle}
\bauthor{\bsnm{Yan}, \binits{Y.}},
\bauthor{\bsnm{Xu}, \binits{Y.}},
\bauthor{\bsnm{Lin}, \binits{Q.}},
\bauthor{\bsnm{Liu}, \binits{W.}},
\bauthor{\bsnm{Yang}, \binits{T.}}:
\batitle{Optimal epoch stochastic gradient descent ascent methods for min-max
  optimization}.
\bjtitle{Advances in Neural Information Processing Systems}
\bvolume{33},
\bfpage{5789}--\blpage{5800}
(\byear{2020})
\end{barticle}
\endbibitem

%%% 24
\bibitem[\protect\citeauthoryear{Du and Hu}{2019}]{du2019linear}
\begin{bchapter}
\bauthor{\bsnm{Du}, \binits{S.S.}},
\bauthor{\bsnm{Hu}, \binits{W.}}:
\bctitle{Linear convergence of the primal-dual gradient method for
  convex-concave saddle point problems without strong convexity}.
In: \bbtitle{The 22nd International Conference on Artificial Intelligence and
  Statistics},
pp. \bfpage{196}--\blpage{205}
(\byear{2019}).
\bcomment{PMLR}
\end{bchapter}
\endbibitem

%%% 25
\bibitem[\protect\citeauthoryear{Shalev-Shwartz and
  Zhang}{2013}]{shalev2013stochastic}
\begin{botherref}
\oauthor{\bsnm{Shalev-Shwartz}, \binits{S.}},
\oauthor{\bsnm{Zhang}, \binits{T.}}:
Stochastic dual coordinate ascent methods for regularized loss minimization.
Journal of Machine Learning Research
\textbf{14}(1)
(2013)
\end{botherref}
\endbibitem

%%% 26
\bibitem[\protect\citeauthoryear{Zhang and Lin}{2015}]{zhang2015stochastic}
\begin{bchapter}
\bauthor{\bsnm{Zhang}, \binits{Y.}},
\bauthor{\bsnm{Lin}, \binits{X.}}:
\bctitle{Stochastic primal-dual coordinate method for regularized empirical
  risk minimization}.
In: \bbtitle{International Conference on Machine Learning},
pp. \bfpage{353}--\blpage{361}
(\byear{2015}).
\bcomment{PMLR}
\end{bchapter}
\endbibitem

%%% 27
\bibitem[\protect\citeauthoryear{Bousquet and
  Elisseeff}{2002}]{bousquet2002stability}
\begin{barticle}
\bauthor{\bsnm{Bousquet}, \binits{O.}},
\bauthor{\bsnm{Elisseeff}, \binits{A.}}:
\batitle{Stability and generalization}.
\bjtitle{The Journal of Machine Learning Research}
\bvolume{2},
\bfpage{499}--\blpage{526}
(\byear{2002})
\end{barticle}
\endbibitem

%%% 28
\bibitem[\protect\citeauthoryear{Shalev-Shwartz
  et~al.}{2010}]{shalev2010learnability}
\begin{barticle}
\bauthor{\bsnm{Shalev-Shwartz}, \binits{S.}},
\bauthor{\bsnm{Shamir}, \binits{O.}},
\bauthor{\bsnm{Srebro}, \binits{N.}},
\bauthor{\bsnm{Sridharan}, \binits{K.}}:
\batitle{Learnability, stability and uniform convergence}.
\bjtitle{The Journal of Machine Learning Research}
\bvolume{11},
\bfpage{2635}--\blpage{2670}
(\byear{2010})
\end{barticle}
\endbibitem

%%% 29
\bibitem[\protect\citeauthoryear{Hardt et~al.}{2016}]{hardt2016train}
\begin{bchapter}
\bauthor{\bsnm{Hardt}, \binits{M.}},
\bauthor{\bsnm{Recht}, \binits{B.}},
\bauthor{\bsnm{Singer}, \binits{Y.}}:
\bctitle{Train faster, generalize better: Stability of stochastic gradient
  descent}.
In: \bbtitle{International Conference on Machine Learning},
pp. \bfpage{1225}--\blpage{1234}
(\byear{2016}).
\bcomment{PMLR}
\end{bchapter}
\endbibitem

%%% 30
\bibitem[\protect\citeauthoryear{Lin et~al.}{2016}]{lin2016generalization}
\begin{bchapter}
\bauthor{\bsnm{Lin}, \binits{J.}},
\bauthor{\bsnm{Camoriano}, \binits{R.}},
\bauthor{\bsnm{Rosasco}, \binits{L.}}:
\bctitle{Generalization properties and implicit regularization for multiple
  passes sgm}.
In: \bbtitle{International Conference on Machine Learning},
pp. \bfpage{2340}--\blpage{2348}
(\byear{2016}).
\bcomment{PMLR}
\end{bchapter}
\endbibitem

%%% 31
\bibitem[\protect\citeauthoryear{Madden et~al.}{2020}]{madden2020high}
\begin{botherref}
\oauthor{\bsnm{Madden}, \binits{L.}},
\oauthor{\bsnm{Dall'Anese}, \binits{E.}},
\oauthor{\bsnm{Becker}, \binits{S.}}:
High-probability convergence bounds for non-convex stochastic gradient descent.
arXiv preprint arXiv:2006.05610
(2020)
\end{botherref}
\endbibitem

%%% 32
\bibitem[\protect\citeauthoryear{Koren and Levy}{2015}]{koren2015fast}
\begin{botherref}
\oauthor{\bsnm{Koren}, \binits{T.}},
\oauthor{\bsnm{Levy}, \binits{K.}}:
Fast rates for exp-concave empirical risk minimization.
Advances in Neural Information Processing Systems
\textbf{28}
(2015)
\end{botherref}
\endbibitem

%%% 33
\bibitem[\protect\citeauthoryear{Zhang et~al.}{2022}]{zhang2022uniform}
\begin{botherref}
\oauthor{\bsnm{Zhang}, \binits{S.}},
\oauthor{\bsnm{Hu}, \binits{Y.}},
\oauthor{\bsnm{Zhang}, \binits{L.}},
\oauthor{\bsnm{He}, \binits{N.}}:
Uniform convergence and generalization for nonconvex stochastic minimax
  problems.
arXiv preprint arXiv:2205.14278
(2022)
\end{botherref}
\endbibitem

%%% 34
\bibitem[\protect\citeauthoryear{Kakade
  et~al.}{2012}]{kakade2012regularization}
\begin{barticle}
\bauthor{\bsnm{Kakade}, \binits{S.M.}},
\bauthor{\bsnm{Shalev-Shwartz}, \binits{S.}},
\bauthor{\bsnm{Tewari}, \binits{A.}}:
\batitle{Regularization techniques for learning with matrices}.
\bjtitle{The Journal of Machine Learning Research}
\bvolume{13}(\bissue{1}),
\bfpage{1865}--\blpage{1890}
(\byear{2012})
\end{barticle}
\endbibitem

%%% 35
\bibitem[\protect\citeauthoryear{Gonen and
  Shalev-Shwartz}{2017}]{gonen2017average}
\begin{barticle}
\bauthor{\bsnm{Gonen}, \binits{A.}},
\bauthor{\bsnm{Shalev-Shwartz}, \binits{S.}}:
\batitle{Average stability is invariant to data preconditioning: Implications
  to exp-concave empirical risk minimization}.
\bjtitle{The Journal of Machine Learning Research}
\bvolume{18}(\bissue{1}),
\bfpage{8245}--\blpage{8257}
(\byear{2017})
\end{barticle}
\endbibitem

%%% 36
\bibitem[\protect\citeauthoryear{Shalev-Shwartz}{2007}]{shalev2007online}
\begin{bbook}
\bauthor{\bsnm{Shalev-Shwartz}, \binits{S.}}:
\bbtitle{Online Learning: Theory, Algorithms, and Applications}.
\bpublisher{Hebrew University}, \blocation{???}
(\byear{2007})
\end{bbook}
\endbibitem

%%% 37
\bibitem[\protect\citeauthoryear{Zhivotovskiy and
  Hanneke}{2018}]{zhivotovskiy2018localization}
\begin{barticle}
\bauthor{\bsnm{Zhivotovskiy}, \binits{N.}},
\bauthor{\bsnm{Hanneke}, \binits{S.}}:
\batitle{Localization of vc classes: Beyond local rademacher complexities}.
\bjtitle{Theoretical Computer Science}
\bvolume{742},
\bfpage{27}--\blpage{49}
(\byear{2018})
\end{barticle}
\endbibitem

%%% 38
\bibitem[\protect\citeauthoryear{Kanade et~al.}{2022}]{kanade2022exponential}
\begin{botherref}
\oauthor{\bsnm{Kanade}, \binits{V.}},
\oauthor{\bsnm{Rebeschini}, \binits{P.}},
\oauthor{\bsnm{Vaskevicius}, \binits{T.}}:
Exponential tail local rademacher complexity risk bounds without the bernstein
  condition.
arXiv preprint arXiv:2202.11461
(2022)
\end{botherref}
\endbibitem

\end{thebibliography}

\end{document}